\documentclass[12pt]{amsart}
\usepackage{graphicx,amssymb,stmaryrd,amsfonts,epsfig,amsthm,a4,amsmath,mathrsfs,url}
\usepackage{vmargin,xypic,amscd}
\usepackage{eufrak}
\usepackage[latin1]{inputenc}
\psfigdriver{dvips}

\theoremstyle{definition}

\theoremstyle{theorem}
\newtheorem{thm}{Theorem}[section]
\newtheorem{cor}[thm]{Corollary}
\newtheorem{lem}[thm]{Lemma}
\newtheorem{prop}[thm]{Proposition}
\newtheorem{prop?}[thm]{Proposition?}
\theoremstyle{definition}
\newtheorem{defn}[thm]{Definition}

\newtheorem{rem}[thm]{Remark}
\newtheorem{que}[thm]{Question}

\newtheorem{exe}[thm]{Example}

\numberwithin{equation}{section}

\newcommand{\Z}{\mathbf{Z}}

\newcommand{\N}{\mathbf{N}}
\newcommand{\R}{\mathbf{R}}
\newcommand{\C}{\mathbf{C}}
\newcommand{\Q}{\mathbf{Q}}
\newcommand{\K}{\mathbf{K}}
\newcommand{\h}{\mathfrak{h}}
\newcommand{\g}{\mathfrak{g}}
\newcommand{\tpr}{\begin{tiny}\noindent Proof:}

\newcommand{\mk}{\mathfrak}

\newcommand{\GL}{\textnormal{GL}}

\newcommand{\SL}{\textnormal{SL}}

\newcommand{\SO}{\textnormal{SO}}

\newcommand{\Cone}{\textnormal{Cone}}

\newcommand{\Exp}{\textnormal{Exp}}
\newcommand{\Rnil}{\textnormal{R}_\textnormal{nil}}
\newcommand{\Rexp}{\textnormal{R}_\textnormal{exp}}

\newcommand{\Expn}{\textnormal{Exp}^\textnormal{nil}}

\newcommand{\Aut}{\textnormal{Aut}}

\begin{document}~



\title{Dimension of asymptotic cones of Lie groups}
\author{Yves de Cornulier}%
\date{\today}

\begin{abstract}We compute the covering dimension the asymptotic cone
of a connected Lie group. For simply connected solvable Lie groups,
this is the codimension of the exponential radical.

As an application of the proof, we give a characterization of connected Lie
groups that quasi-isometrically embed into a non-positively curved
metric space.\end{abstract}
\maketitle

\section{Introduction}

There are many kinds of dimensions associated to a connected Lie
group. The most naive one is the usual dimension; for the point of
view of coarse geometry it is not interesting as it is not a
quasi-isometry invariant, since for instance it would have to vanish for all
compact Lie groups. However, the dimension of $G/K$, where $K$ is a
maximal compact subgroup of $G$, is a quasi-isometry invariant of
$G$ by a theorem of J. Roe \cite[Proposition~3.33 and
Corollary~3.35]{Roe}. Other kinds of dimensions that provide
quasi-isometry invariants are, among others, asymptotic dimension
\cite{Gro,BeDr}, subexponential corank \cite{BuSc1}, hyperbolic
dimension \cite{BuSc2}. Here we focus on the covering dimension of
the asymptotic cone, which was considered in Gromov's book \cite[Chap.~2]{Gro}, in Burillo's paper
\cite{Bur} and in \cite{LePi} in the context of solvable groups, in \cite{BeMi} in the case of the mapping class group, in \cite{DrSm} in connection  with asymptotic Assouad-Nagata dimension, and in \cite{Kl} in the context of non-positively curved manifolds.

The covering dimension $\textnormal{covdim}(X)$ of a topological
space X is the minimal $d$ such that every open covering has a
refinement with multiplicity at most $d+1$. This is a topological
invariant, and as the asymptotic cone (with respect to a given
ultrafilter) is a bilipschitz invariant of quasi-isometry classes of
metric spaces, the cone dimension is a quasi-isometry invariant.

Let us recall the definition of asymptotic cones. For more details
see for instance \cite{Drutu}. Fix a non-principal ultrafilter
$\omega$ on $\N$. Let $X$ be a metric space. Define its asymptotic
precone $\textnormal{Precone}(X,\omega)$ as the following
semi-metric space. As a set, $\textnormal{Precone}(X,\omega)$ is the
set of sequences $(x_n)$ in $X$ such that $(d(x_n,x_0)/n)_{n\ge 1}$
is bounded. The semi-metric is defined as
$d((x_n),(y_n))=\lim_\omega d(x_n,y_n)/n$. The asymptotic cone
$\textnormal{Cone}(X,\omega)$ is the metric space obtained from
$\textnormal{Precone}(X,\omega)$ by identifying points at distance
0. Define the cone dimension $\textnormal{conedim}(X,\omega)$ as the
covering dimension of $\textnormal{Cone}(X,\omega)$.

Besides, recall that, in a simply connected solvable Lie group $G$,
the exponential radical $\Rexp G$ is defined as the set of
exponentially distorted elements. It was introduced by Osin \cite{Osin} 
who
proved that this is a normal closed connected subgroup, and that $G/\Rexp
G$ is the biggest quotient of $G$ with polynomial growth.

\begin{thm}
Let $G$ be a simply connected solvable Lie group. Then
$$\textnormal{conedim}(G,\omega)=\dim(G/\Rexp G)$$
for every non-principal ultrafilter $\omega$ on
$\N$.\label{thm_main}
\end{thm}

The inequality $\le$ follows from work of Osin \cite{Osin} along
with a slight generalization of a result of Burillo \cite{Bur} that
we establish in Section \ref{sec_burillo}.

When the map $G\to G/\Rexp G$ is split, $G/\Rexp G$ embeds
quasi-isometrically into $G$ so that the other inequality is
immediate, as observed in \cite[Proposition~4.4]{LePi}. However in
general this map is not split, even in groups that do not deserve to
be considered as ``pathological". Example
\ref{exe_rexp_non_split_alg} gives a simply connected
group $G$ in which this map is not split, and which moreover is an
algebraic group defined over the rationals, and having an arithmetic
lattice.

\begin{rem}The fact that, for a connected Lie group,
the cone dimension is independent of the ultrafilter $\omega$ is anything
but surprising. However, it is not known if a connected Lie
group can have two non-homeomorphic asymptotic cones; nevertheless
in \cite{KSTT,KrTe} it is proved that an absolutely simple group of real rank at
least two has a unique asymptotic cone up to homeomorphism if and
only if the continuum hypothesis is true.

Nevertheless, there are many connected Lie groups for which the asympotic cone
is known to be unique up to (bilipschitz) homeomorphism. The best-known examples are
groups with polynomial growth, and word hyperbolic Lie groups. There are other examples, such as the group $\textnormal{SOL}(\R)$, for which every asymptotic cone is a ``Diestel-Leader $\R$-graph", namely the ``hypersurface" of equation $b(x)+b(y)=0$ in\footnote{Similarly, a few well-known groups do have the same asymptotic cone as $\textnormal{SOL}(\R)$, such as the solvable Baumslag-Solitar group $\Z[1/n]\rtimes_n\Z$ or the lamplighter group $\Z/n\Z\wr\Z$ ($|n|>1$) . Incidentally, this shows that two groups, one finitely presented and the other not, can have all their asymptotic groups bilipschitz homeomorphic.} the product $T\times T$, where $T$ is a complete universal homogeneous $\R$-tree of degree $2^{\aleph_0}$ and $b$ is a Busemann function on $T$. It is easy to extend such a description to every semidirect product $\R^m\ltimes\R^n$, where the action of $\R^n$ on $\R^m$ is semisimple.

On the other hand it was recently established
\cite[Remark 5.13]{OOS} that there exists a finitely generated
group for which the dimension of the asymptotic cone does depend on the
ultrafilter. More precisely, they construct a finitely generated
(not finitely presented) group having one asymptotic cone homeomorphic to
a $\R$-tree $T$ and therefore of covering dimension one, and one other homeomorphic to
$T\times\R/\Z$ and therefore of covering dimension two.\end{rem}


Before sketching the proof of Theorem \ref{thm_main}, let us give
some corollaries and reformulations.

Although it is well-known that the case of general connected Lie
groups reduces to that of simply connected solvable Lie groups (see
Section \ref{sec_exprad}), it is convenient to have a formula not
appealing to a solvable cocompact group. To write such a formula, we
extend the notion of exponential radical as follows.

First recall \cite[Definition~1.8]{Cor} that a connected Lie group
is called {\em M-decomposed} if every non-compact simple Lie subgroup
centralizes the radical, or equivalently if it is almost the direct
product of a semisimple group with no compact factors and a
solvable-by-compact group, called its amenable radical. If moreover
the amenable radical has polynomial growth, we call the group {\em
P-decomposed}\footnote{It is proved in \cite{CPS} that a connected
Lie group has the {\em Rapid Decay (RD) Property} if and only if it
is P-decomposed.}.

\begin{defn}
If $G$ is a connected Lie group, its exponential radical $\Rexp G$
is the smallest kernel of a map of $G$ onto a P-decomposed connected
Lie group.\label{def_rexp_general}
\end{defn}

In Section \ref{sec_exprad} we check that this definition makes
sense, and extends that of Osin, and we give general results
concerning the exponential radical (see Theorem \ref{thm_rexp}).

\begin{defn}1) If $G$ is a connected Lie group, define its geometric dimension
$\textnormal{geodim}(G)$ as $\dim(G)-\dim(K)$, where $K$ is a
maximal compact subgroup: this is the dimension of a contractible
Riemannian manifold on which $G$ acts properly cocompactly by
isometries.

2) If $S$ is a connected semisimple Lie group with centre $Z$,
define its geometric rank $\textnormal{georank}(S)$ as $r+z$, where
$r$ is the $\R$-rank of $S/Z$, and $z$ is the $\Q$-rank of $Z$ (i.e.
the dimension of the vector space $Z\otimes_\Z\R$).
\end{defn}


\begin{cor}
Let $G$ be a connected Lie group with radical $R$. We have, for
every non-principal ultrafilter $\omega$
$$\textnormal{conedim}(G,\omega)=\textnormal{georank}(G/R)+
\textnormal{geodim}(R/\Rexp G).$$\label{cor_dimcone_lie}
\end{cor}

%




\begin{rem}
It follows from Corollary \ref{cor_dimcone_lie} that if $G$ is a
connected Lie group and $H$ a quotient of $G$, then
$\textnormal{conedim}(G)\ge \textnormal{conedim}(H)$. This is not
true for general finitely generated groups, as we see by taking
$G=F_2$ (free group) and $H=\Z^2$, where $\textnormal{conedim}(G)=1$
and $\textnormal{conedim}(H)=2$. Note also that there exist finitely
generated groups with infinite cone dimension: for instance,
$\Z\wr\Z$ has infinite cone dimension: indeed, it contains
quasi-isometrically embedded free abelian subgroups of arbitrary
large rank.
\end{rem}


The proof of Theorem \ref{thm_main} splits into three steps. Let $G$
denote a simply connected solvable Lie group.

The first paramount step is the case of polynomial growth, and is
already settled: we have to know that
$$\textnormal{conedim}(G/\Rexp G,\omega)=\dim(G/\Rexp G).$$
Indeed, if $H$ is a $d$-dimensional simply connected Lie group of
polynomial growth, then its asymptotic cone is homeomorphic to
$\R^d$: in the essential case where $P$ is nilpotent, this is due to
Pansu \cite{Pan1}, and in general this follows from the nilshadow
construction (see for instance Lemma \ref{lem_polynomial_alg_emb} or
\cite{Bre}). Finally, we have to know that $\R^d$ has covering
dimension $d$ \cite[Theorem~IV.4]{Nag}.



The second step is the inequality
$$\textnormal{conedim}(G,\omega)\le\textnormal{conedim}(G/\Rexp
G,\omega).$$ It relies on Osin's result that the exponential radical
is strictly exponentially distorted, along with a variant of theorem
of Burillo \cite[Theorem~16]{Bur} about the behaviour of the
covering dimension under a map with ultrametric fibers. See the
discussion following Corollary \ref{cor_burillo}.

The last step is the inequality
$$\textnormal{conedim}(G,\omega)\ge\textnormal{conedim}(G/\Rexp
G,\omega).$$ It is obvious when the extension
$$1\to\Rexp G\to G\to G/\Rexp G\to 1$$
is a semidirect product. However, this is not always the case, see
Examples \ref{exe_rexp_non_split_alg} and
\ref{exe_rexp_non_split_non_alg}. In the general case, it relies on
the following result.

\begin{thm}
Let $G$ be a connected solvable Lie group. Fix any non-principal
ultrafilter $\omega$ on $\N$. Then there exists a bilipschitz
embedding $$\textnormal{Cone}(G/\Rexp G,\omega)\to
\textnormal{Cone}(G,\omega).$$ In particular,
$\textnormal{conedim}(G/\Rexp
G,\omega)\le\textnormal{conedim}(G,\omega)$.\label{thm_lift_intro}
\end{thm}

We construct explicitly such an embedding in Section \ref{sec_lift}
when $G/\Rexp G$ is simply connected nilpotent. However in general
$G/\Rexp G$ is only of polynomial growth; nevertheless a {\em
trigshadow} construction, carried out in Section
\ref{sec_trigshadow}, allows to reduce to the case when $G/\Rexp G$
is nilpotent.

Finally, let us give an analogous result for $p$-adic algebraic
groups. If $G$ is a linear algebraic group over an ultrametric local
field of characteristic zero (equivalently, some finite extension of
$\Q_p$), the assumption that $G$ is compactly generated is strong
enough to ensure that the role of the exponential radical is played
by the unipotent radical (see Proposition
\ref{prop_ultrametric_padic}). In this case the Burillo dimension
Theorem \cite[Theorem~16]{Bur} applies directly, and moreover as the
unipotent radical is always split, our lifting construction (Theorem
\ref{thm_lift_intro}) is not needed. See the detailed discussion in
Section \ref{sec_padic}, whose conclusion is the following.

\begin{thm}
Let $G$ be a linear algebraic group over an ultrametric local field
of characteristic zero. Suppose that $G$ is compactly generated and
let $d$ be the dimension of a maximal split torus. Then, for every
non-principal ultrafilter $\omega$,
$$\textnormal{conedim}(G,\omega)=d.$$\label{thm_main_padic}
\end{thm}

Note that it follows from Corollary \ref{cor_dimcone_lie} that if
$G$ is a linear algebraic group over $\R$ or $\C$, with unipotent
radical $U$, and $d$ is the dimension of a maximal split torus, then
$\Rexp G\subset U$ and
$$\textnormal{conedim}(G,\omega)=d+\dim_\R(U/\Rexp G),$$
where $\dim_\R$ means twice the dimension when the ground field is
$\C$.

Let us make a few remarks about other kinds of 
dimension.


\begin{itemize}
\item For every topological space $X$, one can define its (large) inductive 
dimension, 
$\textnormal{Ind}(X)$, which is a non-negative integer (see \cite[Section~I.4]{Nag}). If 
$X$ is 
metrizable (as all asymptotic cones as considered in the paper), then it is known that it coincides with the covering dimension.


\item As the asymptotic cone is defined up to bilipschitz homeomorphism, its 
Hausdorff dimension, which is a non-negative real number, is well-defined.
If $G$ has polynomial growth, it is known to have growth of 
degree $d$ for some (explicit) non-negative integer $d$ \cite{Guiv}, and work of Pansu 
\cite{Pan1} shows that 
the Hausdorff dimension of the asymptotic cone is exactly $d$.

On the other hand, if $G$ is a Lie group or a compactly generated linear 
algebraic group over 
a local field of characteristic zero and if $G$ has non-polynomial growth, 
then by \cite{CT} it contains a 
quasi-isometrically embedded regular trivalent tree, and therefore all its 
asymptotic cone contain a bilispchitz embedded complete universal 
$\R$-tree, everywhere branched of degree $2^{\aleph_0}$. This $\R$-tree 
has the property that every $r$-ball contains uncountably many disjoint 
open $r/2$-balls, and it immediately follows from the definition of 
Hausdorff dimension that the Hausdorff dimension of this $\R$-tree, and 
therefore of all asymptotic cones of $G$, is 
infinite. 
\end{itemize}

Finally, here is an application of our results to embeddings into non-positively curved metric spaces.

Recall that a geodesic metric space is called {\it non-positively curved} or {\it CAT(0)} if for every triple of points $a,b,c$ and every $t\in[0,1]$, the distance of $(1-t)b+tc$ to $a$ is less or equal to the same distance computed from a triangle $a,b,c$ inside the Euclidean plane, with the same edge lengths. For instance, every simply connected Riemannian manifold with non-positive sectional curvature is a CAT(0) metric space.

\begin{thm}
Let $G$ be a simply connected solvable Lie group endowed with a left-invariant Riemannian
metric. Suppose that $G$ is triangulable (i.e. embeds as a subgroup of real upper triangulable matrices, see Lemma \ref{lem_triangulable}). Then the following are equivalent
\begin{itemize}
\item [(1)] $G$ quasi-isometrically embeds into a CAT(0) space;
\item [(2)] $G$ quasi-isometrically embeds into a non-positively curved simply connected symmetric space;
\item [(2')] $G$ bilipschitz embeds into a non-positively curved simply connected symmetric space;
\item [(3)] $G/\Rexp G$ is abelian.
\end{itemize}\label{thm_intro_catzero}
\end{thm}

This is proved in Section \ref{catzero}, where we formulate a statement for all connected Lie groups (although the general case is based, as always, on the case of simply connected Lie groups).

It is in the implication (1)$\Rightarrow$(3) that we make use of Theorem \ref{thm_lift_intro},
combined with the following theorem of Pauls \cite{Pauls}: if $G$ is a non-abelian simply connected nilpotent Lie group, then $G$ has no quasi-isometric embedding into any CAT(0) space.

\begin{rem}
If $G$ is any linear algebraic group over a local field $\K$ of characteristic zero, then any algebraic embedding $G\to\SL_n(\K)$ is quasi-isometric by \cite{Mus} and thus induces a quasi-isometric embedding of $G$ into a CAT(0) space, namely the Euclidean Bruhat-Tits building \cite{BT} of $\SL_n(\K)$.
\end{rem}

\begin{que}
Which non-positively curved simply connected symmetric spaces quasi-isometrically embed into a finite product of trees (resp. binary trees)?
\end{que}

It was proved in \cite{BS} that every word hyperbolic group quasi-isometrically embeds into a finite product of binary trees, and therefore the question has a positive answer for products of rank one non-positively curved simply connected symmetric space. But I do not know the answer for any irreductible non-positively curved simply connected symmetric space of higher rank, e.g. $\SL_3(\R)/\SO_3(\R)$. Of course it makes sense to extend the question to general connected Lie groups.


\section{Trigshadow}\label{sec_trigshadow}

\begin{lem}
Let $G$ be a connected Lie subgroup of $\GL_n(\R)$ and $H$ its
Zariski closure. Then if $G$ has polynomial growth, so does
$H$.\label{lem_poly_zariski}
\end{lem}
\begin{proof}Write $H=DKU$, where $U$ is the unipotent radical, $K$ is
reductive anisotropic, $D$ is an isotropic torus, $[D,K]=1$ and
$K\cap D$ is finite.

Consider the adjoint action of $\GL_n(\R)$ on its Lie algebra. Then
$G$ fixes its Lie algebra $\g$, on which every element of $G$ acts
with eigenvalues of modulus one. Therefore all of $H$ fixes $\g$,
acting with eigenvalues of modulus one. On the other hand, $D$ acts
diagonally with real eigenvalues, therefore acts trivially on $\g$.
It follows that $D$ is, up to a finite central kernel, a direct
factor of $KU$ in $H$, so that $H$ has polynomial growth.
\end{proof}

\begin{lem}
Let $G$ be a connected Lie group with polynomial growth. Then there
exists a connected linear algebraic group $H$ such that
\begin{itemize}\item $H(\R)$ has polynomial growth;
\item $H$ contains no isotropic torus (equivalently: the unipotent radical $U$ of $H$ is
cocompact in $H(\R)$)\item $G$ maps properly (i.e. with compact
kernel) into $H(\R)$, with cocompact, Zariski dense
image.
\item If $K$ is a maximal compact subgroup of $G$,
    then $\dim(U)=\dim(G)-\dim(K)$.
\end{itemize}\label{lem_polynomial_alg_emb}
\end{lem}
The unipotent radical $U$ is known as the {\em nilshadow} of $G$. The construction
given here is close to the original one \cite[Theorem 4.1]{AuGr}.
 
\begin{proof}
By Lemma \ref{lem_solv_lin_rep}, $G$ has a proper mapping into some
$\GL_n(\R)$. Therefore we can suppose that $G$ is a closed subgroup
in $\GL_n(\R)$. Let $B$ be its Zariski closure. Then $B$ also has
polynomial growth by Lemma \ref{lem_poly_zariski}.

Now consider a vector space $V$, viewed as a unipotent algebraic
$\R$-group, with a proper surjective map $G/[G,G]\to V$. Consider
the diagonal mapping $G\to H=B/D\times V$. Then this map is proper.

We have obtained that $G$ maps properly into an algebraic group
$H=KW$ with $K$ compact and $W$ unipotent. View $G$ as a closed
subgroup. Then $GK$ contains $G$ as a cocompact subgroup, and it is
the semidirect product of $K$ by $G\cap W$; in particular the latter
is connected. Accordingly $G\cap W$, being a connected subgroup of
$W$, is an algebraic group. Therefore $G$ is cocompact in the
algebraic $\R$-group $KG$.

To prove the last equality, first observe that we take the quotient
by $(G\cap U)_0$ and therefore suppose that $G\cap U$ is discrete.
Thus $G$ maps with discrete kernel, into $H/U$, which is abelian. So
the Lie algebra of $G$ is abelian, and by connectedness this implies
that $G$ is abelian. Since $G$ is cocompact in $H$, the desired
formula is immediate.
\end{proof}

The following lemma is a particular case of \cite[Lemma 5.2]{Mostow2};
for completeness we include the short proof.

\begin{lem}
Let $G$ be a solvable, connected Lie group $G$ whose derived
subgroup is simply connected (e.g. $G$ itself is simply
connected). Then $G$ has a faithful linear proper (i.e. with
closed image) linear representation.\label{lem_solv_lin_rep}
\end{lem}
\begin{proof}
By Ado's theorem, choose an embedding of $\g$ into
$\mathfrak{gl}_n(\C)$, and deduce a continuous mapping
$G\to\GL_n(\C)$ with discrete kernel $K$. Conjugating if necessary
using Lie-Kolchin's Theorem, we can suppose that $[G,G]$ maps to
the upper unipotent group, in particular $K\cap [G,G]=\{1\}$. Now
choose a faithful proper linear representation of $G/[G,G]$ in
some $\GL_m(\C)$. Then the diagonal map $i:G\to \GL_{n+m}(\C)$ is
an has kernel $K'=K\cap [G,G]$. We claim that it is proper.
Indeed, let $g_n$ be a sequence in $G/K$ such that $i(g_n)\to 1$.
Then $g_n\to 1$ in $G/[G,G]$. Denote $p:G\to G/[G,G]$ the quotient
morphism, and consider a section $j:G/[G,G]\to 1$ (not a
homomorphism) that is continuous at 1. Then
$d_n=g_nj(p(g_n))^{-1}$ is a sequence in $[G,G]$ such that
$i(d_n)\to 1$. On the other hand, as $[G,G]$ is mapped
(conjugating if necessary) to unipotent matrices, it is properly
embedded. Therefore $d_n\to 1$, so that $g_n\to 1$.
\end{proof}


%



\begin{lem}[trigshadow]
Let $G$ be a solvable, connected Lie group $G$ whose derived
subgroup is simply connected. Then $G$ embeds cocompactly into a
group $H$ containing a cocompact connected subgroup $T$, containing
$[G,G]$, that has a faithful proper embedding into some group of
real upper triangular matrices.\label{lem_trigshadow}
\end{lem}
\begin{proof}
Using Lemma \ref{lem_solv_lin_rep}, fix a faithful proper linear
real representation of $G$. Consider its Zariski closure $DKU$,
where $U$ is the unipotent radical, $DK$ a maximal torus, $D$ its
isotropic part and $K$ its anisotropic part. Then set $H=GK$ and
$T_1=H\cap DU$. As $G$ normalizes its Lie algebra $\g$, so does $K$;
therefore $K$ normalizes $G$, so that $H=KG$ is a subgroup; as $K$
is compact and $G$ closed, $H$ is closed. The restriction to $H$ of
the natural map $DKU\to DKU/DU\simeq K/(D\cap K)$ being surjective,
its kernel $T_1$ is cocompact in $H$. Finally, as the Zariski
closure of $T_1$ has no anisotropic torus, it is triangulable. Note
that $[G,G]\subset U$ and therefore $[G,G]\subset T_1$. We claim
that $T$ has finitely many components. Indeed, $H=T_1K$ is connected
and $T_1\cap K=D\cap K$ is finite. We conclude by defining $T$ as
the unit component of $T_1$.
\end{proof}

\begin{lem}
Let $G$ be a connected Lie group. The following are equivalent
\begin{itemize}
 \item[(i)]$G$ embeds into a group of triangular matrices;
 \item[(i')]$G$ embeds properly (i.e. as a closed
 subgroup in the real topology) into a group of triangular matrices;
 \item[(ii)]$G$ embeds into a solvable linear algebraic group with
 no anisotropic torus;
\end{itemize}
Besides, if these conditions are fulfilled, then $G/\Rexp G$ is
nilpotent, i.e. the inclusion of $\Rexp G$ in the stable term of the
descending central series of $G$ is an equality.\label{lem_triangulable}
\end{lem}
\begin{proof}
(i')$\Rightarrow$(i)$\Rightarrow$(ii) is trivial. The implication
(ii)$\Rightarrow$(i) follows from the fact that every split torus is
conjugate to a group of diagonal matrices. Finally to get (i) it
suffices to observe that every connected Lie subgroup of the group
of triangular matrices is closed.

Finally suppose that the conditions are fulfilled: more precisely
suppose (ii) and embed $G$ as a Zariski dense subgroup of a solvable
linear algebraic group $H$ with no anisotropic torus, and with
unipotent radical $U$. Note that the exponential radical $W$ of $G$
is contained in $U$; therefore it is Zariski closed in $H$. The
group $G/W$ is Zariski dense in $H/W$ and has polynomial growth;
therefore $H/W$ has polynomial growth by Lemma
\ref{lem_poly_zariski}. As $H/W$ has no anisotropic torus, it is
therefore the direct product of its unipotent radical by a split
torus. In particular it is nilpotent, hence $G/W$ is nilpotent.
\end{proof}

\begin{defn}
We call a simply connected solvable Lie group $G$ {\it weakly
triangulable} if $G/\Rexp G$ is nilpotent.
\end{defn}

\begin{rem}
Let $\R$ act on $\R^2$ by the one-parameter subgroup
$$\begin{pmatrix}
  e^t\cos t & -e^t\sin t \\
  e^t\sin t & e^t\cos t \\
\end{pmatrix};\;t\in\R.$$
Then $G=\R\ltimes\R^2$ is weakly triangulable but not triangulable.
\end{rem}

\section{Around Burillo's dimension theorem}\label{sec_burillo}

\begin{defn}
A map $f:X\to Y$ between metric spaces has {\em metrically parallel
fibers} if for every $x,y\in X$ there exists $z\in X$ such that
$$d(x,z)=d(f(x),f(z))\text{ and }f(y)=f(z).$$
\end{defn}

\begin{exe}
Let $G$ be a locally compact group generated by a compact subset
$K$, and $G/N$ a quotient of $G$, generated by the image of $K$.
Endow both with the corresponding word metrics. Then the quotient
map $G\to G/N$ has metrically parallel fibers, and so does the
induced map
$$\text{Cone}(G,\omega))\to\text{Cone}(G/N,\omega)),$$
whose fibers are isometric to $\text{Cone}(N,\omega)$, $N$ being
endowed with the restriction of the word metric of $G$.
\end{exe}

\begin{thm}
Let $X\to Y$ be a map between metric spaces, with metrically
parallel fibers. Suppose that all fibers are bilipschitz to
ultrametric spaces, with bilipschitz constants independent of the
fiber. Then $\dim(X)\le\dim(Y)$.
\end{thm}

This theorem is stated in \cite[Theorem~16]{Bur} with the slightly
stronger assumption that all fibers are ultrametric; this is not
sufficient for our purposes. However, the proof follows the same
lines, so we omit some details.

Recall \cite{Nag} that a metric space $Y$ has covering dimension at
most $n$ if and only if there exists a sequence of positive numbers
$(r_k)\to 0$, and a sequence of coverings $(\mathcal{U}_k)$ by open
subsets of diameter at most $r_k$, such that every $x\in Y$ belongs
to at most $n+1$ elements of $\mathcal{U}_k$, and
$\mathcal{U}_{k+1}$ is a refinement of $\mathcal{U}_k$ for all $k$.

\begin{proof}
We can suppose that there exists $\lambda>0$ and, for every fiber
$X_y=f^{-1}(\{y\})$, a map from the fiber to an ultrametric space
$\Upsilon_y$ satisfying, for all $x,x'\in X_y$
$$d(x,x)\le d(f(x),f(x'))\le \lambda d(x,x').$$

Suppose that $Y$ has covering dimension at most $n$. Consider a
sequence $(\mathcal{U}_k)$ of open coverings of $Y$ as above.
Extracting if necessary, we can suppose that $r_{k+1}\le
r_k/(3\lambda+1)$ for all $k$.

Write $\mathcal{U}_k=(U_k^\alpha)_\alpha$. For every $\alpha$ pick
$y_\alpha\in U_k^\alpha$. As $\Upsilon_{y_\alpha}$ is ultrametric,
it can be covered by pairwise disjoint open balls
$B_k^{\alpha\beta}$ of radius $3r_k$. More precisely, for all $x,x'$
in $\Upsilon_{y_\alpha}$, either they are in the same ball and
$d(x,x')<3r_k$, either they are in two distinct balls and
$d(x,x')\ge 3r_k$.

Consider the corresponding covering ${B'_k}^{\alpha\beta}$ of
$F_{y_\alpha}$ by disjoint open subsets. For all $x,x'$ in
$F_{y_\alpha}$, either they belong to the same subset of the
covering and $d(x,x')< 3\lambda r_k$, or they belong to two distinct
subsets and $d(x,x')\ge 3r_k$.

Define $V_{k}^{\alpha\beta}$ as the set of elements in
$f^{-1}(U_\alpha)$ that are at distance $\le r_k$ from
${B'_k}^{\alpha\beta}$. This is, $\alpha$ being fixed, a family of
open subsets of $X$; any two elements in the same subset are at
distance $<(3\lambda+2)r_k$, and any two elements in two distinct
subsets in this family ($\alpha$ being fixed) are at distance $\ge
r_k$, in particular for $\alpha$ fixed they are pairwise disjoint.
This implies that any element of $X$ belongs to at most $n+1$
elements of the family ($\alpha$ being no longer fixed).

Moreover, this is a covering of $X$, denoted $\mathcal{V}_k$:
indeed, if $x\in X$, then $x\in f^{-1}(U_\alpha)$ for some $\alpha$,
and therefore there exists $x'\in F_\alpha$ such that
$d(x,x')=d(f(x),y_\alpha)\le r_k$, so $x'$ belongs to some
${B'_k}^{\alpha\beta}$, implying $x\in V_{k}^{\alpha\beta}$.

Finally it remains to show that $\mathcal{V}_{k+1}$ is a refinement
of $\mathcal{U}_{k}$. If $V_{k+1}^{\alpha\beta}$ belongs to
$\mathcal{V}_{k+1}$, then $U_{k+1}^\alpha$ is contained in a element
$U_k^\gamma$ of $\mathcal{U}_k$. Now, we have just shown that if
$x\in V_{k+1}^{\alpha\beta}$, then $x\in V_k^{\gamma\delta_x}$ for
some $\delta_x$. Now we claim that all $\delta_x$ are equal to one
single element $\delta$, proving $V_{k+1}^{\alpha\beta}\subset
V_k^{\gamma\delta}$. Indeed, if $\delta_x\neq\delta_{x'}$, then $x$
and $x'$ are at distance at least $r_k$; on the other hand they are
at distance $<(3\lambda+2)r_{k+1}$, contradicting the assumption
$r_{k+1}\le r_k/(3\lambda+2)$.
\end{proof}

We call a metric space $X$ quasi-ultrametric if it satisfies a
quasi-ultrametric inequality: for some $C<\infty$ and for all
$x,y,z\in X$,
$$d(x,z)\le \max(d(x,y),d(y,z))+C.$$

Note that this immediately implies that any asymptotic cone of $X$
is ultrametric. Important examples are log-metrics, that is, when
$(X,d)$ is a metric space, the new metric space $(X,\log(1+d))$.

\begin{cor}
Let $G$ be a locally compact compactly generated group, and $N$ a
closed normal subgroup. Suppose that $N$, endowed with the word
metric of $G$, is quasi-isometric to a quasi-ultrametric space.
Then, for every non-principal ultrafilter $\omega$,
$$\textnormal{conedim}(G,\omega)\le\textnormal{conedim}(G/N,\omega).$$\label{cor_burillo}
\end{cor}

Let $G$ be a simply connected solvable Lie group. By Osin's Theorem,
the restriction $d$ of the word metric (or Riemannian metric) to the
exponential radical is equivalent to the metric $d'$ defined by the
length $\log(|\cdot|_{R_{\exp}(G)}+1)$.

Corollary \ref{cor_burillo} then implies
\begin{equation}\dim(\text{Cone}(G,\omega))\le
\dim(\text{Cone}(G/\Rexp G,\omega)).
\label{eq_dimcone2}\end{equation}



\section{Splittings of the exponential radical}


As we mentioned in the introduction, it is obvious that if the
exponential radical is split, i.e. when the extension
$$1\to\Rexp G\to G\to G/\Rexp G\to 1$$
is a semidirect product, then the converse of Inequality
(\ref{eq_dimcone2}) holds.

However it is not true in general that the exponential radical is
split. Here is an example, which is moreover $\Q$-algebraic with no
$\Q$-split torus and which therefore contains an arithmetic lattice.

\begin{exe}
Let $H_3(R)$ denote the Heisenberg group over the ring $R$. If $K$
is a field of characteristic zero, the automorphism group of the
algebraic group $H_3(K)$ is easily seen (working with the Lie
algebra) to be $K$-isomorphic to $\GL_2(K)\ltimes K^2$ (where the
$K^2$ corresponds to inner automorphisms). If an automorphism
corresponds to $(A,v)\in\GL_2(K)\ltimes K^2$, its action on the
center of $H_3(K)$ is given by multiplication by $\det(A)^2$.

Think now at $K=\R$ and pick a one-parameter diagonalizable subgroup
$(\alpha_t)$ inside $\SL_2(\R)$, viewed as a subgroup of
$\text{Aut}H_3(\R)$. If we want to fix a $\Q$-structure, choose
$\alpha_1=\begin{pmatrix}
  2 & 1 \\
  1 & 1 \\
\end{pmatrix}$ so that this group is $\Q$-anisotropic (although
$\R$-isotropic). This defines a semidirect product $\R^*\ltimes
H_3(\R)$. Define $G_1=(\R^*\ltimes H_3(\R))\times H_3(\R)$. Its
center is the product of both centers and is therefore isomorphic to
$\R\times\R$. Let $Z$ be the diagonal of $\R\times\R$, and define
$$G=G_1/Z.$$ Then $\Rexp G$ is the left-hand copy of $H_3(\R)$. The
quotient $G/\Rexp G$ is isomorphic to $R^3$. However it is not
split: indeed otherwise, $G$ would have $\R^2$ as a direct factor,
while the centre of $G$ is isomorphic to $\R$. Note that $G$ is
triangulable.\label{exe_rexp_non_split_alg}
\end{exe}

Obviously, if $G$ is a connected solvable Lie group with
$\dim(G/\Rexp G)\le 1$ (note that $\dim(G/\Rexp G)=0$ only if $G$ is
compact), then the exponential radical is split.

Also, if $G$ is a connected solvable Lie group that is algebraic
over $\R$ (that is, the unit component in the real topology of
$H(\R)$ for some linear algebraic $\R$-group $H$), and if
$\dim(G/\Rexp G)\le 2$, then it is not hard to see that the
exponential radical is split. Here is a non-algebraic example where
$\dim(G/\Rexp G)=2$ and $\Rexp G$ is not split.

\begin{exe}
Keep the one parameter subgroup in $\Aut(H_3)$ of Example
\ref{exe_rexp_non_split_alg}. Denote by $H'_3$ another copy of
$H_3$, and fix any homomorphism of $H'_3$ onto $\R$. Thus define a
group $G_1=H'_3\ltimes H_3$, where $H'_3$ acts on $H_3$ through the
homomorphism onto $\R$ and the one-parameter subgroup fixed above.
(Note that $G_1$ is definitely non-algebraic because $H'_3$ plays
both a semisimple and a unipotent role.) Let $Z$ be the diagonal of
the center $\R\times\R$ of $G_1$, and define $G=G_1/Z$. Then $\Rexp
G=H_3$ has codimension 2 in $G$ and is not split; note that $G$ is
trangulable.\label{exe_rexp_non_split_non_alg}
\end{exe}


\section{Liftings in weakly triangulable groups}\label{sec_lift}



The purpose of this section is to prove that Equality
(\ref{eq_dimcone2}) holds without assuming that the exponential
radical is split. This is a consequence of the following theorem.

\begin{thm}
Let $G$ be a connected solvable Lie group. Fix any non-principal
ultrafilter $\omega$ on $\N$. Then there exists a bilipschitz
embedding $$\textnormal{Cone}(G/\Rexp G,\omega)\to
\textnormal{Cone}(G,\omega).$$ In particular,
$\textnormal{conedim}(G/\Rexp
G,\omega)\le\textnormal{conedim}(G,\omega)$.
\end{thm}


\begin{proof}

The second statement follows from the first as if $Y$ is a closed
subset of a topological space $X$, then obviously
$\textnormal{covdim}(Y)\le\textnormal{covdim}(X)$.



By Lemma \ref{lem_trigshadow}, we can suppose that the intersection
$N$ of the descending central series of $G$ coincides with the
exponential radical. Denote by $\g$, $\mk{n}$ the corresponding Lie
algebras. It is known \cite[Chapitre~7, p.19-20]{Bou} that there
exists a ``Cartan subalgebra", which is nilpotent subalgebra $\h$ of
$\g$ such that $\h+\mk{n}=\g$. (Examples
\ref{exe_rexp_non_split_alg} and \ref{exe_rexp_non_split_non_alg}
show that we cannot always demand in addition $\h\cap\mk{n}=0$.) Set
$\mk{w}=\h\cap\mk{n}$ and choose a complement subspace $\mk{v}$ of
$\mk{w}$ in $\mk{n}$. Denote by $G,H,N,W$ the corresponding simply
connected Lie groups, and by $\nu$ the restriction to $\mk{v}$ of
the projection $\g\to\g/\mk{n}$.

Define a map $\psi:G/N(=H/W)\to H\subset G$ as
$$\psi=\exp_{H}\circ\nu^{-1}\circ{\exp_{H/W}}^{-1}.$$
This obviously a continuous section of the projection $G\to G/N$,
but it can be checked that this is not a homomorphism nor is large
Lipschitz in general. Nevertheless, we are going to show that $\psi$
induces a map $\tilde{\psi}:\Cone(G/N)\to\Cone(G)$. Obviously
$\tilde{\psi}$ is expansive (i.e. increases distances), and in
particular it is injective. Eventually we are going to show that
$\tilde{\psi}$ is bilipschitz.



To avoid lengthy notation, it is convenient to identify every
(simply connected) nilpotent Lie group with its Lie algebra. From
that point of view, $\psi$ is simply the inclusion of $\mk{v}$ in
$\g$, where $\g$ is endowed with the word metric $|\cdot|_G$ of $G$
and $\mk{v}$ is endowed with the word metric $|\cdot|_{G/N}$ of
$G/N$, through the identification of $\mk{g}/\mk{n}$ and $\mk{v}$ by
$\nu$. Write $\mk{h}=\mk{v}\oplus\mk{w}$. Write $f\preceq g$ if
there exists constants $\alpha,\beta>0$ such that $f\le \alpha
g+\beta$. Fix a norm $\|\cdot\|$ on $\g$.

\begin{lem}
If $(x,w)\in\mk{h}=\mk{v}\oplus\mk{w}$, then
$$|(x,w)|_G\preceq |x|_{G/H}+\log(1+\|w\|)$$\label{lem_long}
\end{lem}
\begin{proof}
Let $|\cdot|$ denote the word length on $G$, $|\cdot|_N$ that on
$N$, and let $\pi$ denote the projection $G\to G/N$. Fix a compact
symmetric generating subset $S$ of $H$. If $g\in H$, write $\pi(g)$
as an element of minimal length with respect to $\pi(S)$:
$\pi(g)=\pi(s_1)\dots\pi(s_m)$. Set $h=s_1\dots s_m$. As $h^{-1}g$
belongs to the exponential radical, the inequality $|g|\le
|h|+|h^{-1}g|$ yields, in view of \cite[Theorem 1.1(3)]{Osin},
$$|g|\preceq |\pi(g)|_{G/N}+\log(1+|h^{-1}g|_N)\quad(g\in H)$$
and therefore
$$|g|\preceq |\pi(g)|_{G/N}+\log(1+\|h^{-1}g\|).$$
As $\h$ is nilpotent, its group law is given by a polynomial and we
have an upper bound $\|h^{-1}g\|\le A(\|g\|+1)^k(\|h\|+1)^k$ for
some constants $A,k\ge 1$. This implies
$$|g|\preceq |\pi(g)|_{G/N}+\log(1+\|g\|)+\log(1+\|h\|).$$
As $H$ is nilpotent, we have $\|h\|\le A'|h|_H^{k'}$ for some
constants $A',k'\ge 1$ (this follows from \cite[Lemma II.1]{Guiv}).
Now $|h|_H=m\le |g|_H\preceq \|g\|$. Therefore
$$|g|\preceq |\pi(g)|+\log(1+\|g\|).$$

Take now $(x,u)$ as in the statement of the lemma. Then the above
formula yields
$$|(x,u)|\preceq |x|+\log(1+\|(x,u)\|),$$ so that
$$|(x,u)|\preceq |x|+\log(1+\|x\|)+\log(1+\|u\|).$$
As $\log(1+\|x\|)\preceq |x|_H\preceq |x|$, we get
$$|(x,u)|\preceq |x|+\log(1+\|u\|).\qedhere$$
\end{proof}



Lemma \ref{lem_long} implies that $\psi$ preserves the word length,
and therefore induces a map
$\tilde{\psi}:\textnormal{Precone}(G/N)\to\textnormal{Cone}(G)$.

On the other hand, observe that $\psi$ is a section of the
projection $G\to G/N$, which decreases distances; in particular
$\psi$, and therefore $\tilde{\psi}$, increases distances. However
let us show that $\tilde{\psi}$ is bilipschitz.


Denote by $\cdot$ the group laws in $\mk{v}$ and $\mk{h}$ coming
from $H/W$ and $H$; although they do not coincide in restriction to
$\mk{v}$ there will be no ambiguity as elements of $\mk{h}$ will
always be written as pairs $(x,w)\in\mk{v}\oplus\mk{w}$. Note that
in both spaces, the inverse law is simply given by $x^{-1}=-x$ (this
simply translates the fact that $\exp(-x)=\exp(x)^{-1}$).

Write the group law in $\mk{h}$, restricted to $\mk{v}$, as
$$(x,0)\cdot(y,0)=(x\cdot y,P(x,y)),$$
where $P$ is a polynomial vector-valued function.


Then $$(x,0)^{-1}\cdot (y,0)=(x^{-1}\cdot y,P(-x,y)),\quad\forall
x,y\in\mk{v}.$$ Introduce the following notation:

If $a\ge 1$ and $b\ge 0$, write $u\le_{a,b} v$ if $u\le av+b$. Then
by Lemma \ref{lem_long}
$$|(x,0)^{-1}\cdot (y,0)|_{G}\le_{a,b}|x^{-1}\cdot y|_{G/N}+\log(1+\|P(-x,y)\|),\quad\forall
x,y\in\mk{v}.
$$
Write $P=\sum P_iv_i$ where $P_i$'s are real-valued polynomials and
$v_i$ are unit vectors. Then
$$|(x,0)^{-1}\cdot (y,0)|_{G}\le_{a,b}|x^{-1}\cdot y|_{G/N}+\sum\log(1+|P_i(-x,y)|),\quad\forall
x,y\in\mk{v}.$$

Now let us work in the asymptotic cone. All limits and cones are
with respect to a fixed ultrafilter $\omega$. Any linearly bounded
sequences $(x_n)$ and $(y_n)$ in $\mk{v}$ (i.e.
$|x_n|_{H/W},|y_n|_{H/W}=O(n)$) define elements $\underline{x}$ and
$\underline{y}$ in both $\Cone(H/W)$ and $\Cone(G)$. We have

\begin{eqnarray}
\nonumber
d(\underline{x},\underline{y})_{\Cone(G)}&=&\lim\frac{1}{n}|(x_n,0)^{-1}(y_n,0)|_{G}\\
\nonumber &\le_{a,0} & \lim\frac{1}{n}|{x_n}^{-1}\cdot
y_n|_{G/N}+\lim\frac{1}{n}\sum\log(1+|P_i(-x_n,y_n)|)\\
\nonumber &=&
d(\underline{x},\underline{y})_{\Cone(H/W)}+\sum\lim\frac{1}{n}\log(1+|P_i(-x_n,y_n)|).
\end{eqnarray}

Now, as $(x_n)$ and $(y_n)$ are linearly bounded,
$|P_i(-x_n,y_n)|=O(n^d)$ where $d$ is the the total degree of $P$.
Thus $\lim\frac{1}{n}\log(1+|P_i(-x_n,y_n)|)=0$ for all $i$, and we
obtain

$$ d(\underline{x},\underline{y})_{\Cone(G)}\le_{a,0}
d(\underline{x},\underline{y})_{\Cone(H/W)}.$$ This shows that
$\tilde{\psi}$ is lipschitz, and therefore is bilipschitz. In
particular, it factors through a map
$$\textnormal{Cone}(G/N)\to\textnormal{Cone}(G),$$
ending the proof.\end{proof}


\begin{rem}We have shown that the asymptotic cone of every simply
connected solvable Lie group $G$ contains a bilipschitz-embedded
copy of $\Cone(G/\Rexp G)$. It is natural to ask if $G$ itself
contains a quasi-isometrically embedded copy of $G/\Rexp(G)$. For
instance, if $G$ is the group of Example
\ref{exe_rexp_non_split_non_alg}, does $G$ contain a flat?
\end{rem}



\section{Exponential radical}\label{sec_exprad}

The exponential radical was introduced by Osin in \cite{Osin} for
simply connected solvable Lie groups. We extend here the definition
to all connected Lie groups. Recall that an element $g$ of a locally
compact, compactly generated group $G$ is {\em exponentially
distorted} if $$|g^n|=O(\log n),$$ where $|\cdot|$ denotes word
length. Let $\Exp(G)$ denote the set of exponentially distorted
elements (which includes torsion elements).

\begin{thm}[Osin]
Let $G$ be a simply connected, solvable Lie group. Then
\begin{enumerate}
    \item $\Exp(G)$ is a closed, connected subgroup of
$G$ (necessarily characteristic), called the exponential radical of
$G$ and denoted $R_{\exp}(G)$.
    \item \label{item_osin_logmetric}The exponential radical is
    {\em strictly exponentially distorted}, that is the
restriction of the metric $|\cdot|$ to $R_{\exp}(G)$ is equivalent
to the log-metric $\log(|\cdot|_{R_{\exp}(G)}+1)$, where
$|\cdot|_{R_{\exp}(G)}$ denotes the intrinsic word metric of
$R_{\exp}(G)$.
\item $G/R_{\exp}(G)$ is the biggest quotient of $G$ with polynomial
growth.\end{enumerate}\label{thm_osin}
\end{thm}


\begin{rem} This definition of exponential radical is not suitable
for general connected Lie groups, in view of the two following
examples. In $\SO(2)\ltimes\R^2$, which is solvable but not simply
connected, and in $\SL_2(\C)$, which is simply connected but not
solvable, the subset $\Exp(G)$ is neither a subgroup, nor is
closed.\end{rem}

\begin{defn}
Let $G$ be a connected Lie group, and $R$ its (usual) radical.
Define its exponential radical $\Rexp(G)$ so that $G/\Rexp(G)$ is
the biggest quotient of $G$ that is locally isomorphic to a direct
product of a semisimple group and a group with polynomial growth.
\end{defn}

It is easily checked that $\Rexp(G)$ is the closed normal subgroup
generated by $[S_{nc},R]$, where $S_{nc}$ is a maximal semisimple
subgroup with no compact factors. (Actually $[S_{nc},R]$ is a normal
subgroup, see Lemma \ref{lem_SncR_exp_tordu}.)

The following lemma shows that this extends Osin's definition.

\begin{lem}
Let $G$ be a simply connected, solvable Lie group. Then
$G/R_{\exp}(G)$ is the biggest quotient of $G$ with polynomial
growth.\label{lem_osin_rexp_polyq}
\end{lem}
\begin{proof}
Let $M$ be the smallest connected normal subgroup of $G$ such that
$G/M$ has polynomial growth. As $G/M$ is a simply connected solvable
Lie group with polynomial growth, it satisfies $\Exp(G/M)=1$. This
implies that the image of $\Exp(G)=\Rexp(G)$ in $G/M$ is trivial,
that is, $\Rexp(G)\subset M$.

By \cite[Theorem 1.1(4)]{Osin}, $\Rexp(G/\Rexp(G))=1$, and therefore
by \cite[Proposition 3.2]{Osin}, $G/\Rexp(G)$ has polynomial growth,
so that $\Rexp(G)=M$.

Now suppose that $N$ is a closed, normal subgroup of $G$ such that
$G/N$ has polynomial growth. Let $N_0$ be the connected component of
$G$. Then $G/N_0$ is a central extension of $G/N$ which has
polynomial growth, so that $G/N_0$ also has polynomial growth. By
definition of $M$ this implies $M\subset N_0$ and therefore
$M\subset N$.
\end{proof}

\begin{thm}
Let $G$ be a connected Lie group,
\begin{itemize}
    \item[(1)] $\Rexp G$ is a closed connected characteristic subgroup of
    $G$.
    \item[(2)] $\Rexp G$ is contained in the nilpotent radical of
    $G$.
    \item[(3)]\label{item_descr_rexp}  Denote $R$ the radical of $G$ and $S_{nc}$ a
    maximal semisimple subgroup with no compact factors. Then
        $$\Rexp G=\overline{\Rexp(R)[S_{nc},R]}.$$
    \item[(4)] Suppose that $G$ is algebraic, that $U$ is its
    unipotent radical, and $L_{nc}$ a maximal reductive isotropic
    subgroup. Then
        $$\Rexp G=[L_{nc},U].$$
    \item[(5)] $\Rexp G$ is strictly exponentially distorted in $G$.
\end{itemize}\label{thm_rexp}
\end{thm}

We need a series of lemmas.

\begin{lem}Let $A$ be a connected abelian Lie
group, $K$ its maximal torus, and $F$ be a closed subgroup of $A$
such that $A\cap K=1$. Then $F$ is contained in a direct factor of
$K$ in $A$.\label{lem_direct_factor}
\end{lem}
\begin{proof}
Taking the quotient by the connected component $F_0$, we can suppose
that $F$ is discrete, and thus is free abelian. Write
$A=\tilde{A}/\Gamma$, and lift $F$ to a subgroup $\tilde{F}$ of
$\tilde{A}$; we have $\tilde{F}\cap\Gamma=\{0\}$. As
$\tilde{F}+\Gamma$ is closed and countable, it is discrete. This
implies that $\tilde{F}\cap\text{Vect}_\R(\Gamma)=\{0\}$. Thus
$\tilde{F}$ is contained in a complement subspace of
$\text{Vect}_\R(\Gamma)$ and we are done.
\end{proof}

The following lemma is well-known but we found no reference.

\begin{lem}
Let $K$ be a field of characteristic zero, $\g$ a finite dimensional
Lie algebra, and $\mk{n}$ its nilpotent radical. Then $\g/\mk{n}$ is
reductive, i.e. is the direct product of a semisimple Lie subalgebra
and an abelian one. \label{lem_qt_rnil_reductive}
\end{lem}
\begin{proof}
Fix an embedding of $\g$ into $\mk{gl}_n$. Let $H\subset\GL_n$ be
the connected normalizer of $\g$, $R$ the radical of $H$ and $U$ its
unipotent radical. Then $[H,R]\subset U$. Denoting by Gothic letters
the corresponding Lie algebras, we have
$[\mk{h},\mk{r}]\subset\mk{u}$. Denoting by $\mk{r}'$ the radical of
$\g$, we have $\mk{r}'\subset\mk{r}$ and therefore
$[\mk{g},\mk{r}']\subset\mk{u}\cap\g$. So $[\mk{g},\mk{r}']$ is a
nilpotent ideal in $\g$ and is thus contained in its nilpotent
radical.
\end{proof}

\begin{defn}
If $G$ is a group and $H$ a subgroup, we say that $H$ is cocentral
in $G$ if $G$ is generated by $H$ and its centre
$Z(G)$.\label{def_cocentral}
\end{defn}

Note that if $H$ is cocentral in $G$, every normal subgroup of $H$
remains normal in $G$.

The following lemma is quite standard.

\begin{lem}
Let $G$ be a connected Lie group and $N$ its nilradical. Then $G$
has a cocompact subgroup $L$ (which can be chosen normal if $G$ is
solvable), containing $N$, and embedding cocompactly cocentrally
into a connected solvable Lie group $M$ whose fundamental group
coincides\footnote{If $G$ is a topological group and $H$ a subgroup,
we say that $\pi_1(G)$ coincides with $\pi_1(H)$ if the natural map
$\pi_1(H)\to\pi_1(G)$ is an isomorphism.} with that of $N$. In
particular, if $G$ has no nontrivial compact normal torus, then $M$
is simply connected. Moreover, $M$ decomposes as a direct product
$Q\times\R^n$ for some $n$ and some closed subgroup $Q$ of $L$; in
particular $\Rexp(M)\subset G$.\label{lem_solvable_cocompact}
\end{lem}
\begin{proof}Working inside $G/N$ and taking the preimage in $G$
allows us to work in $H=G/N$. Let $R$ be the radical of $H$; it is
central in $H$ by Lemma \ref{lem_qt_rnil_reductive}. Let $S$ be the
semisimple Levi factor in $H$, and $\bar{S}$ its closure. The centre
of the semisimple part of $G$ has a subgroup of finite index whose
action by conjugation on $N$ is trivial. Thus if $\bar{S}\cap
R=\overline{Z(S)}\cap R$ contained a torus, the preimage of this
torus in $G$ would be nilpotent and this would contradict the
definition of $N$. Thus $\bar{S}\cap R$ has a subgroup of finite
index whose intersection with the maximal torus $K$ of $R$ is
trivial, which is therefore contained in a direct factor $L$ of $K$
in $R$ by Lemma \ref{lem_direct_factor}. We obtain that
$H_1=L\bar{S}$ is a connected, cocompact subgroup of $H$ whose
radical $R_1$ is simply connected.

The group $H_1/Z(H_1)$ is semisimple with trivial centre and thus
has a connected, solvable, simply connected cocompact subgroup $B$.
Let $C$ be the preimage of $B$ in $H_1/R_1$. Then the covering
$H_1/R_1\to H_1/Z(H_1)$ induces a covering $C\to B$, with fibre
$Z(H_1)/R\simeq Z(H_1/R)$. As $B$ is simply connected, we deduce
that $C$ is isomorphic to $B\times Z(H_1/R_1)$; in particular we
identify $B$ with the unit component of $C$. The group $Z(H_1/R_1)$
has a free abelian subgroup $\Lambda$ of finite index, which embeds
as a lattice in some vector space $V$. Let $Q$ be the preimage of
$B\subset H_1/R_1$ in $G$. Then $Q$ is a closed, connected, solvable
subgroup of $G$, whose fundamental group coincides with that of $N$.

On the other hand, let $\tilde{S}$ be a Levi factor of $G$. Then the
quotient map $G\to G/R$ maps $Z(\tilde{S})$ onto $Z(H_1/R_1)$. Now
$Z(\tilde{S})$ contains a subgroup of finite index which is central
in $G$; let $\Upsilon$ be its closure. Then $Q\Upsilon$ has finite
index in the preimage of $C$ in $G$; in particular it is a closed,
cocompact solvable subgroup. As $Q$ contains $N$, it also contains
$\Upsilon_0$; so there exists a discrete, free abelian subgroup
$\Gamma$ of $\Upsilon$ such that $\Gamma\cap Q=1$ and $Q\Gamma$ has
finite index in $Q\Upsilon$. Now $L=Q\times\Gamma$ embeds
cocompactly cocentrally in $M=Q\times(\Gamma\otimes_\Z\R)$.
\end{proof}

\begin{defn}If $G$ is a connected Lie group, define
$\Expn G$ as $\Exp G\cap \Rnil G$.
\end{defn}

\begin{lem}
Let $G$ be a connected Lie group $R$ its radical, $S$ a Levi factor,
and $S_{nc}$ its non-compact part. Then $[S_{nc},R]$ is normal in
$G$, and is contained in $\Expn(G)$.\label{lem_SncR_exp_tordu}
\end{lem}
\begin{proof}

Let us first check that $[S_{nc},R]$ is normal. First observe that
$[S_{nc},R]=[S_{nc},[S_{nc},R]]$: this is obtained by noting that
$S_{nc}=[S_{nc},S_{nc}]$ and using the Jacobi identity. Then
$$[G,[S_{nc},R]]=[G,[S_{nc},[S_{nc},R]]]
\subset[S_{nc},[[S_{nc},R],G]]+[[S_{nc},R],[G,S_{nc}]]$$
$$\subset [S_{nc},R]+[[S_{nc},R],[S_{nc},R]]\subset [S_{nc},R].$$

Let $N$ be the nilpotent radical. It follows from Lemma
\ref{lem_qt_rnil_reductive} that $[S_{nc},R]$ is nilpotent, and
therefore contained in $N$.

Using Lemma \ref{lem_solvable_cocompact}, consider a cocompact
overgroup $M$ of a cocompact subgroup $L$ of $G$ containing $N$,
with $M$ solvable simply connected, and $\Rexp M\subset G$. Set
$H=\Rexp M\cap N$. It is a closed subgroup of $N$, stable under
taking square roots, and therefore connected. Let $V$ be the
centralizer of $N$ in $G$ modulo $H$ (that is, the preimage in $G$
of the centralizer of $N/H$ in $G/H$). Then $V$ is normal in $G$,
and $\Rexp M\subset V$. We deduce that $M/V$ has polynomial growth.
Therefore $G/V$ has polynomial growth. This implies that $S_{nc}$ is
contained in $V$. In particular $[S_{nc},N]\subset H$. It follows
that

$$[S_{nc},R]=[S_{nc},[S_{nc},R]]\subset [S_{nc},N]\subset H.$$

\end{proof}

\begin{lem}
$\Expn G$ is a closed characteristic subgroup of $G$, containing
$\Rexp G$.\label{lem_expn}
\end{lem}

\begin{proof} Let us prove that $\Expn G$ is a closed subgroup (obviously
characteristic). Let $K$ be the maximal compact subgroup of $\Rnil
G$. By Lemma \ref{lem_solvable_cocompact}, there exists a cocompact
subgroup $L$ of $G/K$, containing $\Rnil(G/K)$, embedding
cocompactly into a simply connected solvable Lie group $M$. Note
that as $G/K$ and $M$ have a common cocompact subgroup,
$\Exp(G/K)\cap L=\Exp(M)\cap L=\Exp M$.  Since $L\supset
\Rnil(G/K)$, we deduce that $\Expn(G/K)$ coincides with $\Exp(M)\cap
\Rnil(G/K)$. It is therefore, by Theorem \ref{thm_osin}, a closed
subgroup of $\Rnil(G/K)$. Now $K$ is central in $G$ and therefore
$\Rnil(G/K)=\Rnil(G)/K$; it follows that that the preimage of
$\Expn(G/K)$ in $G$ coincides with $\Expn(G)$; this is therefore a
closed subgroup.

Now let us prove that $\Rexp G\subset \Expn G$. To show this, using
(3) we only have to check that $\Rexp R\subset \Expn G$ and
$[S_{nc},R]\subset \Expn G$.

The second inclusion is Lemma \ref{lem_SncR_exp_tordu}.

For the first inclusion, we have to show that $R/\Expn G$ has
polynomial growth. As $\Expn R\subset\Expn G$, it suffices to show
that $R/\Expn R$ has polynomial growth.

So let us prove that for any connected solvable Lie group $G$,
$G/\Expn G$ has polynomial growth. We begin as in the beginning of
this proof. The main new feature is that now $L$ is normal in $G$
(see Lemma \ref{lem_solvable_cocompact}). Therefore
$\Expn(G/K)=\Exp(M)$ is normal in $G$. As $M/\Exp M$ has polynomial
growth, so does $(G/K)/\Expn(G/K)=G/K\Expn(G)$.
\end{proof}

\begin{proof}[Proof of Theorem \ref{thm_rexp}]~
\begin{itemize}
    \item[(1)] This is obvious; connectedness is due to the fact that
P-decomposability only depends on the Lie algebra.
    \item[(2)] As the quotient of $G$ by its nilpotent radical $\Rnil G$ is
    a reductive Lie group (see Lemma \ref{lem_qt_rnil_reductive}), $\Rexp G$ is contained in $\Rnil G$.
    \item[(3)] By definition $\Rexp G$ is obviously the closure of
    the normal subgroup generated by $\Rexp(R)$ and $[S_{nc},R]$. So
    we only have to check that $\overline{\Rexp(R)[S_{nc},R]}$ is
    normal in $G$, and this is a consequence of the fact that
    $[S_{nc},R]$ is normal (Lemma \ref{lem_SncR_exp_tordu}).
    \item[(4)] Noting that $[L_{nc},U]$ is closed, using
    (3) allows to reduce to the case when $G$ is
    solvable. Now for $G$ solvable, $\Rexp G$ is the kernel of the
    biggest quotient with polynomial growth. Therefore, denoting by
    $D$ a maximal split torus ($D=L_{nc}$), it follows from the Lie algebra
    characterisation of type R that $[D,G]$, and hence $[D,U]$,
    is contained in $\Rexp G$. Conversely, if $[D,U]=1$, then it
    immediately follows from a Levi decomposition that $G$ is
    nilpotent-by-compact. Therefore all what remains to prove is
    that $[D,U]$ is a normal subgroup. This is checked by means of
    the Lie algebra. First, using that the action of $D$ on
    $\mk{u}$ is semisimple, one checks that $[D,[D,U]]=[D,U]$. Then,
    using this along with the Jacobi Identity, one obtains that $[D,U]$
    is normalized by $U$. Thus it is normal.
    \item[(5)]\label{item_expn_str_dist}Use the notation as in the
    beginning of the proof of Lemma \ref{lem_expn}. We have
    inclusions $\Rexp G/K\subset\Rexp M\subset M$. By Osin's Theorem
    (Theorem \ref{thm_osin}), $\Rexp M$ is
    strictly exponentially distorted in $M$. On the other hand, as
    both $\Rexp G/K$ and $\Rexp M$ are nilpotent, we have polynomial
    distortion: for some $a\ge 1$, we have, for $g\in\Rexp G/K$,
    $$|g|_{\Rexp G/K}^{1/a}\preceq |g|_{\Rexp M}\preceq |g|_{\Rexp
    G/K}.$$ This implies that $\Rexp G/K$ is strictly exponentially
    distorted in $G$, and therefore in $G/K$. Thus $\Rexp G$
    is strictly exponentially
    distorted in $G$.\qedhere
\end{itemize}
\end{proof}

\begin{lem}Let $G$ be a connected Lie group, and $N$ a closed normal
solvable subgroup. Then
$$\textnormal{conedim}(G,\omega)\ge
\textnormal{conedim}(G/N,\omega).$$\label{lem_dimcone_qt_resol}
\end{lem}
\begin{proof}By Lemma \ref{lem_solvable_cocompact}, $G/N$ has a
cocompact solvable subgroup $H/N$ embedding cocompactly into a
connected solvable Lie group $L/N$. Therefore we can suppose $G$
solvable.

We can suppose that $G$ and $G/N$ have no non-trivial compact normal
subgroup. Now taking a cocompact simply connected normal subgroup in
$G/N$ and taking the unit component of its preimage in $G$, we can
suppose that $G$ and $G/N$ are both simply connected solvable Lie
groups. Then the result follows from Theorem \ref{thm_main}.
\end{proof}

\begin{proof}[Proof of Corollary \ref{cor_dimcone_lie}]
Set $c(G)=\textnormal{conedim}(G,\omega)$ and
$\rho(G)=\textnormal{georank}(G/R)+ \textnormal{geodim}(R/\Rexp G)$.

By Corollary \ref{cor_burillo} and using Theorem \ref{thm_rexp}(5),
we have $c(G/\Rexp G)\le c(G)$; by Lemma \ref{lem_dimcone_qt_resol}
this is an equality. On the other hand, it is clear that
$\rho(G)=\rho(G/\Rexp G)$. Thus to prove that $\rho=c$ we are
reduced to the case when $\Rexp G=1$.

Thus suppose $\Rexp G=1$, i.e. $G$ is P-decomposed. Let $M$ denote
the amenable radical of $G$, namely the subgroup generated by $R$
and compact simple factors. Set $S=G/S$, and $Z$ a finite index
subgroup of its centre, free abelian of rank $n$. Then $S/Z$ has a
cocompact simply connected solvable subgroup, whose preimage in $S$
is denoted by $B$. As $B/Z$ is simply connected, $B$ is the direct
product of its unit component $B_0$ and $Z$. Moreover, as $G$ is
M-decomposed, the preimage of $Z$ in $G$ is central. So lift $Z$ to
a subgroup $Z'$ of $G$. Let $H$ be the preimage of $B$ in $G$. Then
$H=Z'\times H_0$. It follows from Theorem \ref{thm_main} that for
connected Lie groups, $c$ is additive under direct products. Thus
$c(H)=n+c(H_0)$. Now as $G$ is M-decomposed and $B_0$ has trivial
centre, $H_0$ is isomorphic to the direct product of $S/Z$ and $M$.
Thus $c(G)=n+c(S/Z)+c(M)$. On the other hand, it is immediate that
$\rho(G)=n+\rho(S/Z)+\rho(M)$. Thus we are reduced to prove $c=\rho$
in two cases: semisimple groups with finite centre and groups with
polynomial growth.

Suppose that $G$ is semisimple with finite centre. Then $G$ has a
cocompact simply connected solvable group $AN$, where $D$ is
isomorphic to $\R^r$, where $r=\R\textnormal{-rank}(G)$, and acts
faithfully on the Lie algebra of $N$ by diagonal matrices and no
fixed points. In particular, $N$ is the exponential radical of $DN$,
and by Theorem \ref{thm_main},
$c(DN)=\R\textnormal{-rank}(G)=\rho(G)$.

Suppose that $G$ has polynomial growth. Then by lemma
\ref{lem_polynomial_alg_emb}, $G$ is quasi-isometric to a simply
connected nilpotent Lie group of dimension
$d=\textnormal{geodim}(G)=\rho(G)$, and the asymptotic cone of $G$
is homeomorphic to $\R^d$ by Pansu's work \cite{Pan1}.
\end{proof}

\section{The $p$-adic case}\label{sec_padic}

In this section we deal with groups $G$ of the form
$\mathbf{G}(\Q_p)$ where $\mathbf{G}$ is a connected linear
algebraic group defined over $\Q_p$ for some prime $p$. For the sake
of simplicity, we call $G$ a connected linear algebraic
$\Q_p$-group.

To define the exponential radical, the notion of exponentially
distorted elements is not relevant. Indeed, it is not hard to show
that in a connected linear algebraic $\Q_p$-group $G$, every element
is either non-distorted or elliptic, and all unipotent elements are
elliptic.

Fix a connected linear algebraic $\Q_p$-group $G$, and take a Levi
decomposition $G=H\ltimes U$ where $U$ is the unipotent radical, and
$H$ is a maximal reductive subgroup \cite{Mostow}. Decompose $H$ as
an almost product $T=SK$, where $S$ is contains all isotropic
factors (simple or abelian) and $K$ contains all anisotropic
factors.

Recall that $G$ is compactly generated if and only if $U=[S,U]$
\cite[Th\'eor\`eme~13.4]{BorelTits}. In this case, the role of the
exponential radical is played by the unipotent radical, in view of
the following proposition, which is probably folklore.

\begin{prop}
Let $G$ be a compactly generated connected linear algebraic group
defined over $\Q_p$, and $U$ its unipotent radical. Then the
restriction to $U$ of the word metric of $G$ is equivalent to an
ultrametric. More precisely, if $G$ is solvable with no non-trivial
anisotropic torus, there exists a left-invariant metric on $G$ that
is equivalent to the word metric, and which is ultrametric in
restriction to $U$.\label{prop_ultrametric_padic}
\end{prop}
\begin{proof}
If we replace $G$ by a cocompact solvable subgroup containing $U$,
then the unipotent radical of the new group contains $U$. Therefore
we can suppose that $G$ is solvable and without anisotropic torus.
It follows that $G$ has a faithful representation into a group of
upper triangular matrices in $\Q_p$. Now S.~Mustapha proved that an
algebraic embedding between compactly generated linear algebraic
groups is always quasi-isometric. Therefore we are reduced to the
case of the upper triangular group $G$ of size $n\times n$ over
$\Q_p$, where $U$ is the set of upper unipotent matrices.

Fix the field norm $|\cdot|$ on $\Q_p$, defined as $|x|=e^{-v_p(x)}$
where $v_p$ is the $p$-valuation. Let us define a real-valued map
$\|\cdot\|$ on $U$ as follows: if $A=(a_{ij})\in U$, set
$\|A\|=\max\{|a_{ij}|^{1/(j-i)}:1\le i<j\le n$. Then $\|\cdot\|$ is
an ultralength on $U$, that is satisfies the inequality
$\|AB\|\le\max(\|A\|,\|B\|)$ and the equalities $\|1\|=0$ and
$\|A^{-1}\|=A$. We omit the details.


Now fix any function $\ell:\R_+\to\R_+$ satisfying $\ell(0)=0$ and

For $x\ge 0$, set $\ell(x)=\max(0,\log(x))$; it satisfies the
inequality $\ell(xy)\le\ell(x)+\ell(y)$.

Denote by $D$ the group of diagonal matrices, so that $G=D\ltimes
U$. Set, for $d=\textnormal{diag}(d_1,\dots,d_n)\in D$,
$\|d\|=\max(|d_1|,\dots,|d_n|)$. Note that for all $(d,u)\in D\times
U$, we have $\|d^{-1}ud\|\le\|d\|^2\|u\|$.



Define, for $g=du$, where $(d,u)\in D\times U$,
$|g|=\max(2\ell(\|d\|),\ell(\|u\|)$. Then $|\cdot|$ is subadditive:
indeed
$$|d_1u_1d_2u_2|=|d_1d_2.d_2^{-1}u_1d_2u_2|=\max(2\ell(\|d_1d_2\|),\ell(\|d_2^{-1}u_1d_2u_2\|))$$
$$\le \max(2\ell(\|d_1\|)+2\ell(\|d_2\|),\ell(\|d_2^{-1}u_1d_2\|),\ell(\|u_2\|))$$
$$\le \max(2\ell(\|d_1\|)+2\ell(\|d_2\|),2\ell(\|d_2\|)+\ell(\|u_1\|),\ell(\|u_2\|))\le |d_1u_1|+|d_2u_2|$$

As $|\cdot|$ is bounded on compact subsets, it follows that it is
dominated by the word length. It remains to check that the word
length is dominated by $|\cdot|$. On $D$, this is obvious, so it
suffices to check this on $U$. Suppose that $u\in U$ satisfies
$|u|\le n$. Let $w\in D$ be the matrix
$\textnormal{diag}(1,p,p^2,\dots,p^{n-1})$. Then $\|w^{-n}uw\|\le
1$. As $\{u\in U:\|u\|\le 1\}$ is compact, this shows that $|\cdot|$
is equivalent to the word length on $G$. Now $|\cdot|$ is
ultrametric on $U$ and this completes the proof.
\end{proof}

Let $G$ be a compactly generated linear algebraic group over a local
field of characteristic zero as above, and $U$ its unipotent
radical. Let $d$ be the dimension of a maximal split torus in $G$.
Then $G$ has a cocompact subgroup $H$ containing $U$ that is
solvable with no non-trivial anisotropic torus, and whose maximal
split tori have dimension $d$. By Proposition
\ref{prop_ultrametric_padic}, $H$ has a left-invariant metric
equivalent to the word metric, for which the restriction to $U$ (and
hence to every coset of $U$)is ultrametric. Therefore
\cite[Theorem~16]{Bur} applies directly to obtain

\begin{equation}\dim(\text{Cone}(H,\omega))\le\dim(\text{Cone}(H/U,\omega)).
\label{eq_dimconepadic}\end{equation}

Of course, we can also use Corollary \ref{cor_burillo} to get it.


Now $H$ has a Levi decomposition: $H=D\ltimes U$ where $D$ is a
split torus, isomorphic to $H/U$. Thus the inclusion of $D$ in $H$
is a quasi-isometric embedding, and therefore induces a bilipschitz
embedding of $\text{Cone}(H/U,\omega)$ into $\text{Cone}(H,\omega)$.

Thus (\ref{eq_dimconepadic}) is an equality. Moreover, $D$ is
isomorphic to ${\K^*}^d$, and $\K^*$ is isomorphic to the direct
product of $\Z$ and of the 1-sphere in $\K$, which is compact. Thus
the asymptotic cone of $D$ is homeomorphic to $\R^d$, and thus we
get, recalling that $H$ is cocompact in $G$

$$\dim(\text{Cone}(G,\omega))=d.$$

\section{Embeddings into non-positively curved metric spaces}\label{catzero}

Define a connected Lie group $G$ as quasi-abelian if it acts properly cocompactly on a Euclidean space (equivalently, if it has polynomial growth and its nilshadow is abelian). Define $G$ as quasi-reductive if its amenable radical $A$ is quasi-abelian and is, up to local isomorphism, a direct factor (note that this is stronger than P-decomposability, see just before Definition \ref{def_rexp_general}).

\begin{thm}
Let $G$ be a connected Lie group endowed with a left-invariant Riemannian
metric. Then the following are equivalent
\begin{itemize}
\item [(1)] $G$ quasi-isometrically embeds into a CAT(0) space;
\item [(2)] $G$ quasi-isometrically embeds into a non-positively curved simply connected symmetric space;
\item [(3)] $G/\Rexp G$ is quasi-reductive.
\end{itemize}
If moreover $G$ is simply connected and solvable, these are also equivalent to:
\begin{itemize}
\item [(2')] $G$ bilipschitz embeds into a non-positively curved simply connected symmetric space;
\end{itemize}\label{thm_emb_catzero}
\end{thm}

In the essential case when $G$ is simply connected and triangulable, (3) is equivalent to: ``$G/\Rexp G$ is abelian" and we get the statement of Theorem \ref{thm_intro_catzero}.

The following lemma may be of independent interest.

\begin{lem}
Let $G$ be a simply connected solvable Lie group. Consider any embedding $i:G\to S$, with closed image, into a solvable Lie group $S$. Denote by $p$ the quotient morphism $G\to G/\Rexp G$.Then the product morphism $i\times p:G\to S\times G/\Rexp G$ is a bilipschitz embedding (where groups are endowed with left-invariant Riemannian metrics).\label{prod_embed}
\end{lem}
\begin{proof}

As a group morphism, $i\times p$ is Lipschitz; as it is injective it is locally bilipschitz; so it suffices to prove that $i\times p$ is a quasi-isometric embedding. Endow our groups with their word length relative to suitable compact generating subsets.


Suppose that $|i\times p(g)|\le n$. Then $|p(g)|\le n$. Thus $p(g)=p(s_1)\dots p(s_n)$ with $s_i\in W$, our compact generating subset of $G$. Set $h=g(s_1\dots s_n)^{-1}$. Then $|g|\le |h|+n$, $h\in\Rexp G$, and $|i(h)|\le 2n$.

Now $|h|\le C_1\log(\|\log(h)\|)+C_1$ by Lemma \ref{estim}. The Lie algebra map induced by $i$ is linear injective, therefore bilipschitz. Therefore, $|h|\le C_1\log(\|\log(i(h))\|)+C_2$. It follows, using Lemma \ref{estim} again, that $|h|\le C_3|i(h)|+C_3\le C_3(2n+1)$. Thus $|g|\le (2C_3+1)n+C_3$ and we are done.
\end{proof}

\begin{lem}
Let $G$ be a simply connected solvable Lie group and $N$ its exponential radical. Take $H\in\mk{n}$, the Lie algebra of $N$, and $h=\exp(H)$. Then we have, for constants independent of $H$
$$C\log(\|H\|+1)-C\le |h|\le C'\log(\|H\|+1)+C'.$$\label{estim}
\end{lem}
\begin{proof}
Let $|\cdot|_N$ denote the intrinsic word length in $N$. Then by Osin's Theorem (Theorem \ref{thm_osin}(\ref{item_osin_logmetric})),
$$C_1\log(|h|_N+1)-C_1\le |h|\le C'_1\log(|h|_N+1)+C'_1.$$
On the other hand, as a consequence of \cite[Lemma II.1]{Guiv} we have,
$$C_2\|H\|^{1/d}-C_2\le |h|_N\le C'_2\|H\|+C'_2\qedhere$$
\end{proof}


\begin{proof}[Proof of Theorem \ref{thm_emb_catzero}]
(2')$\Rightarrow$(2)$\Rightarrow$(1) is trivial.

First suppose $G$ simply connected solvable.
Then (3)$\Rightarrow$(2') follows from Lemma
\ref{prod_embed}, where we choose $S$ as the upper triangular matrices in some $\SL_n(\C)$.

Finally if (1) is satisfied, then the asymptotic cone of $G$ embeds bilipschitz into a CAT(0) space. But the asymptotic cone of $G$ contains by Theorem \ref{thm_lift_intro} a bilipschitz embedded copy of the asymptotic cone of $G/\Rexp G$, which embeds quasi-isometrically into a CAT(0) space only if the nilshadow of $G/\Rexp G$ is abelian by Pauls' Theorem \cite{Pauls}.

The general case follows from the following observation: let $G$ be any connected Lie group. Write $G/\Rexp G=SA$, where $A$ is its amenable radical and $S$ is semisimple without compact factors. Then $G/\Rexp G$ has a cocompact closed subgroup of the form $R\times\Z^n\times A$, where $R$ is a cocompact subgroup of $S/Z(S)$ that is simply connected and solvable. As $R/\Rexp R$ is abelian, we see that $G/\Rexp G$ is quasi-reductive if and only if $A$ is quasi-abelian.

Therefore, if $G$ satisfies (3), $A$ is quasi-abelian and accordingly $G$ is quasi-isometric to $H=R\times\Z^n\times B/K$, where $B,K$ are normal closed subgroups of $A$ with $B$ solvable cocompact, $K$ compact, and $B/K$ simply connected. Now $H$ is solvable simply connected and $H/\Rexp H$ is quasi-abelian, so $H$ and hence $G$ satisfies (2).

Conversely if $G$ satisfies (1), then so does $H$, so by the solvable simply connected case, $H/\Rexp H$ is quasi-abelian. This means that $A$ is quasi-abelian, and therefore $G/\Rexp G$ is quasi-reductive.
\end{proof}

\medskip

\noindent{\bf Acknowledgments.} I thank E.~Breuillard, A.~Eskin, P.~Pansu, and
especially R.~Tessera for valuable discussions.

\end{document}